\documentclass[11pt,a4paper]{article}
\usepackage{graphicx}
\usepackage{amsfonts}
\usepackage{latexsym}

\newcommand{\beql}[1]{\begin{equation}\label{#1}}
\newcommand{\eeq}{\end{equation}}
\newcommand{\comment}[1]{}

\newcommand{\eqref}[1]{{\rm (\ref{#1})}}

\newcommand{\Abs}[1]{{\left|{#1}\right|}}

\newcommand{\Norm}[1]{{\left\|{#1}\right\|}}

\newcommand{\Qed}{\ \\\mbox{$\Box$}}

\newcommand{\Set}[1]{{\left\{{#1}\right\}}}

\newcommand{\RR}{{\mathbb R}}
\newcommand{\CC}{{\mathbb C}}
\newcommand{\ZZ}{{\mathbb Z}}

\newcommand{\one}{{\bf 1}}

\newcommand{\inner}[2]{{\langle #1, #2 \rangle}}
\newcommand{\Inner}[2]{{\left\langle #1, #2 \right\rangle}}
\newcommand{\dens}{{\rm dens\,}}

\newcommand{\ft}[1]{\widehat{#1}}

\newcounter{open}
\newcounter{dfn}
\def\thedfn{\arabic{dfn}}
\newenvironment{dfn}{
  \sf
  \vskip 0.10in
  \refstepcounter{dfn}
  \noindent{\bf Definition \thedfn \ }
}{\vskip 0.10in}

\newcounter{obs}
\def\theobs{\arabic{obs}}

\newcounter{thm}
\newcounter{mysec}
\newcounter{mysubsec}[mysec]
\newtheorem{theorem}{Theorem}

\begin{document}

\begin{center}
{\Large\bf A class of non-convex polytopes that admit no\\
	orthonormal basis of exponentials}\\
\ \\
{\sc Mihail N. Kolountzakis} and {\sc Michael Papadimitrakis}\\
\ \\
\small January 2001\\
\end{center}

\begin{abstract}
A conjecture of Fuglede states that a bounded measurable set $\Omega\subset\RR^d$,
of measure $1$, can tile $\RR^d$
by translations if and only if the Hilbert space
$L^2(\Omega)$ has an orthonormal basis consisting of exponentials
$e_\lambda(x) = \exp 2\pi i\inner{\lambda}{x}$.
If $\Omega$ has the latter property it is called {\em spectral}.
Let $\Omega$ be a polytope in $\RR^d$ with the following property:
there is a direction $\xi \in S^{d-1}$ such that, of all the polytope
faces perpendicular to $\xi$, the total area of the faces pointing in
the positive $\xi$ direction is more than the total area of the faces
pointing in the negative $\xi$ direction.
It is almost obvious that such a polytope $\Omega$ cannot tile space
by translation.
We prove in this paper that such a domain is also not spectral,
which agrees with Fuglede's conjecture.
\end{abstract}

Let $\Omega$ be a measurable subset of $\RR^d$, which
we take for convenience to be of measure $1$.
Let also $\Lambda$ be a discrete subset of $\RR^d$.
We write
\begin{eqnarray*}
e_\lambda(x) &=& \exp{2\pi i \inner{\lambda}{x}},\ \ \ (\lambda,x\in\RR^d),\\
E_\Lambda &=& \Set{e_\lambda:\ \lambda\in\Lambda} \subset L^2(\Omega).
\end{eqnarray*}
The inner product and norm on $L^2(\Omega)$ are
$$
\inner{f}{g}_\Omega = \int_\Omega f \overline{g},
\ \mbox{ and }\ 
\Norm{f}_\Omega^2 = \int_\Omega \Abs{f}^2.
$$

\begin{dfn}
The pair $(\Omega, \Lambda)$ is called a {\em spectral pair}
if $E_\Lambda$ is an orthonormal basis for $L^2(\Omega)$.
A set $\Omega$ will be called {\em spectral} if there is
$\Lambda\subset\RR^d$ such that 
$(\Omega, \Lambda)$ is a spectral pair.
The set $\Lambda$ is then called a {\em spectrum} of $\Omega$.
\end{dfn}

\noindent
{\bf Example:} If $Q_d = (-1/2, 1/2)^d$ is the cube of
unit volume in $\RR^d$ then
$(Q_d, \ZZ^d)$ is a spectral pair
($d$-dimensional Fourier series).

We write $B_R(x) = \Set{y\in\RR^d:\ \Abs{x-y}<R}$.
\begin{dfn} (Density)\\
The discrete set $\Lambda \subset \RR^d$ has {\em density} $\rho$, and we write
$\rho = \dens \Lambda$, if we have
$$
\rho = \lim_{R\to\infty} {\#(\Lambda\cap B_R(x)) \over \Abs{B_R(x)}},
$$
uniformly for all $x\in\RR^d$.
\end{dfn}

We define translational tiling for complex-valued functions below.
\begin{dfn}
Let $f:\RR^d\to\CC$ be measurable and $\Lambda \subset \RR^d$ be a discrete set.
We say that {\em $f$ tiles with $\Lambda$ at level $w\in\CC$}, and sometimes
write ``$f + \Lambda = w \RR^d$'', if
\beql{tiling}
\sum_{\lambda\in\Lambda} f(x - \lambda) = w,
  \ \ \mbox{for almost every (Lebesgue) $x\in\RR^d$},
\eeq
with the sum above converging absolutely a.e.
If $\Omega \subset \RR^d$ is measurable we say that {\em $\Omega + \Lambda$ is a tiling}
when $\one_\Omega + \Lambda = w\RR^d$, for some $w$.
If $w$ is not mentioned it is understood to be equal to $1$.
\end{dfn}

\noindent
{\bf Remark 1}\\
If $f \in L^1(\RR^d)$, $f\ge0$, and $f+\Lambda = w\RR^d$, then
the set $\Lambda$ has density
$$
\dens\Lambda = {w\over\int f}.
$$

The following conjecture is still unresolved in all dimensions
and in both directions.

\noindent
{\bf Conjecture:} (Fuglede \cite{F74}) If $\Omega \subset \RR^d$
is bounded and has Lebesgue measure $1$ then $L^2(\Omega)$ has
an orthonormal basis of exponentials if and only if
there exists $\Lambda \subset \RR^d$ such that $\Omega + \Lambda = \RR^d$ is a tiling.

Fuglede's conjecture has been confirmed in several cases.
\begin{enumerate}
\item
Fuglede \cite{F74} shows that if $\Omega$ tiles with $\Lambda$ being
a lattice then it is spectral with the dual lattice $\Lambda^*$
being a spectrum.
Conversely, if $\Omega$ has a lattice $\Lambda$ as a spectrum
then it tiles by the dual lattice $\Lambda^*$.
\item
If $\Omega$ is a convex non-symmetric domain (bounded, open set) then,
as the first author of the present paper has proved \cite{K00},
it cannot be spectral. It has long been known that convex domains
which tile by translation must be symmetric.
\item
When $\Omega$ is a smooth convex domain it is clear that it admits
no translational tilings.
Iosevich, Katz and Tao \cite{IKT} have shown that
it is also not spectral.
\item
There has also been significant progress in dimension $1$ (the
conjecture is still open there as well) by \L aba \cite{La,Lb}.
For example, the conjecture has been proved in dimension $1$ if
the domain $\Omega$ is the union of two intervals.
\end{enumerate}

In this paper we describe a wide class of, generally non-convex, polytopes
for which Fuglede's conjecture holds.
\begin{theorem}\label{th:main}
Suppose $\Omega$ is a polytope in $\RR^d$ with the following property:
there is a direction $\xi \in S^{d-1}$ such that
$$
\sum_{i} \sigma^*(\Omega_i) \neq 0.
$$
The finite sum is extended over all faces $\Omega_i$ of $\Omega$
which are orthogonal to $\xi$ and
$\sigma^*(\Omega_i) = \pm \sigma(\Omega_i)$,
where $\sigma(\Omega_i)$ is the surface measure of $\Omega_i$
and the $\pm$ sign depends upon whether the outward unit normal
vector to $\Omega_i$ is in the same or opposite direction with $\xi$.

Then $\Omega$ is not spectral.
\end{theorem}
Such polytopes cannot tile space by translation for the following,
intuitively clear, reason.
In any conceivable such tiling the set of positive-looking
faces perpendicular to $\xi$ must be countered by an equal
area of negatively-looking $\xi$-faces, which is impossible
because there is more (say) area of the former than the latter.

It has been observed in recent work on this problem
(see e.g.\ \cite{K00}) that
a domain (of volume $1$) is spectral with spectrum $\Lambda$
if and only if $\Abs{\ft{\chi_\Omega}}^2 + \Lambda$ is a
tiling of Euclidean space at level $1$.
By Remark 1 this implies that $\Lambda$ has density $1$.

By the orthogonality of $e_\lambda$ and $e_\mu$
for any two different $\lambda$ and $\mu$ in $\Lambda$,
it follows that
\beql{eq1}
\ft{\chi_\Omega}(\lambda-\mu) = 0.
\eeq
It is only this property, and the fact that any spectrum
of $\Omega$ must have density $1$, that are used in the
proof.

\noindent
{\bf Proof of Theorem \ref{th:main}.}

\noindent
The quantities $P, Q, N, \ell$ and $K$, which are introduced in
the proof below, will depend only on the domain $\Omega$.
(The letter $K$ will denote several different constants.)

Suppose that $\Lambda$ is a spectrum of $\Omega$.
Define the Fourier transform of $\chi_\Omega$ as
$$
\ft{\chi_\Omega}(\eta) = \int_\Omega e^{-2\pi i \inner{x}{\eta}}\,dx.
$$
By an easy application of the divergence theorem we get
$$
\ft{\chi_\Omega}(\eta) = -{1\over i\Abs{\eta}}
	\int_{\partial\Omega} e^{-2\pi i \inner{x}{\eta}}
	 \Inner{{\eta\over\Abs{\eta}}}{\nu(x)}\,d\sigma(x),
\ \ \ \eta\neq 0,
$$
where $\nu(x) = (\nu_1(x),\ldots,\nu_d(x))$ is the outward unit normal
vector to $\partial\Omega$ at $x \in \partial\Omega$ and
$d\sigma$ is the surface measure on $\partial\Omega$.

From the last formula we easily see that for some $K\ge1$
\beql{eq2}
\Abs{\nabla\ft{\chi_\Omega}(\eta)} \le {K \over \Abs{\eta}},
\ \ \ \Abs{\eta}\ge 1.
\eeq
Without loss of generality we assume that $\xi=(0,\ldots,0,1)$.
Hence
$$
\ft{\chi_\Omega}(t \xi) = -{1\over it} \int_{\partial\Omega}
	e^{-2\pi i t x_d}\nu_d(x)\,d\sigma(x).
$$
Now it is easy to see that each face of the polytope other
than any of the $\Omega_i$s contributes $O(t^{-2})$
to $\ft{\chi_\Omega}(t\xi)$ as $t\to\infty$.
Therefore
\beql{eq3}
\Abs{\ft{\chi_\Omega}(t\xi) +
 {1\over it}\sum_i e^{-2\pi i \lambda_i t} \sigma^*(\Omega_i)}
	\le {K \over t^2},\ \ \ t\ge 1,
\eeq
where $\lambda_i$ is the value of $x_d$
for $x=(x_1,\ldots,x_d)\in\Omega_i$.

Now define
$$
f(t) = \sum_i \sigma^*(\Omega_i) e^{-2\pi i \lambda_i t},\ \ \ t\in\RR.
$$
$f$ is a finite trigonometric sum and has the following properties:
\begin{itemize}
\item[(i)]
$f$ is an almost-periodic function.
\item[(ii)]
$f(0) \neq 0$ by assumption. Without loss of generality assume $f(0) = 1$.
\item[(iii)]
$\Abs{f'(t)} \le K$, for every $t\in\RR$.
\end{itemize}
By (i), for every $\epsilon>0$ there exists an $\ell>0$ such that
every interval of $\RR$ of length $\ell$ contains a translation number
$\tau$ of $f$ belonging to $\epsilon$:
\beql{eq4}
\sup_t \Abs{f(t+\tau) - f(t)} \le \epsilon
\eeq
(see \cite{B32}).

Fix $\epsilon>0$ to be determined later
($\epsilon=1/6$ will do) and the corresponding $\ell$.
Fix the partition of $\RR$ in consecutive intervals of length $\ell$,
one of them being $[0,\ell]$.
Divide each of these $\ell$-intervals into $N$ consecutive
equal intervals of length $\ell/N$,
where
$$
N > {6 K \ell \sqrt{d-1} \over \epsilon}.
$$

In each $\ell$-interval there is at least one ${\ell\over N}$-interval
containing a number $\tau$ satisfying \eqref{eq4}.
For example, in $[0, \ell]$ we may take $\tau = 0$ and
the corresponding ${\ell\over N}$-interval to be $[0, \ell/N]$.

Define the set $L$ to be the union of all these ${\ell\over N}$-intervals
in $\RR$.
Then $L\xi$ is a copy of $L$ on the $x_d$-axis.
Construct $M$ by translating copies of the cube $[0,\ell/N]^d$
along the $x_d$-axis
so that they have their $x_d$-edges on the ${\ell\over N}$-intervals of
$L\xi$.

The point now is that there can be no two $\lambda$s of $\Lambda$
in the same translate of $M$, at distance $D > {2K\over\epsilon}$
from each other.
Suppose, on the contrary, that
$$
\lambda_1, \lambda_2 \in \Lambda,\ \ \Abs{\lambda_1-\lambda_2}\ge D,
\ \ \lambda_1, \lambda_2 \in M+\eta.
$$
Then $\lambda_1 = t_1\xi + \eta + \eta_1$,
$\lambda_2 = t_2\xi + \eta + \eta_2$,
for some $t_1, t_2 \in L$, $\eta_1, \eta_2 \in \RR^d$ with
$$
\Abs{\eta_1}, \Abs{\eta_2} < {\ell\over N}\sqrt{d-1} < {\epsilon\over6K}.
$$
Hence, $\lambda_1-\lambda_2 = (t_1-t_2)\xi + \eta_1-\eta_2$
and an application of the mean value theorem together
with \eqref{eq1} and \eqref{eq2} gives
$$
\Abs{\ft{\chi_\Omega}((t_1-t_2)\xi)} \le 
 {3K\over\Abs{t_1-t_2}}\Abs{\eta_1-\eta_2}.
$$
From \eqref{eq3} we get
$$
\Abs{f(t_1-t_2)} \le 3K\Abs{\eta_1-\eta_2} + {K \over \Abs{t_1-t_2}}
  < 2\epsilon.
$$
Now, since $t_1, t_2 \in L$, there exist $\tau_1, \tau_2$ satisfying
\eqref{eq4} so that
$$
\Abs{\tau_1-t_1}, \Abs{\tau_2-t_2} < {\ell\over N}
$$
and hence (by (iii))
$$
\Abs{f(\tau_1-\tau_2)-f(\tau_1-t_2)},
\Abs{f(\tau_1-t_2)-f(t_1-t_2)} < K{\ell\over N} < \epsilon.
$$
Therefore
\begin{eqnarray*}
2\epsilon & > & \Abs{f(t_1-t_2)} \\
 &\ge& \Abs{f(0)} - \Abs{f(0)-f(-\tau_2)}
                  -\Abs{f(-\tau_2)-f(\tau_1-\tau_2)}\\
 &\ &\ \ \ -\Abs{f(\tau_1-\tau_2)-f(\tau_1-t_2)}
          -\Abs{f(\tau_1-t_2) - f(t_1-t_2)}\\
 &\ge& 1 -\epsilon -\epsilon -\epsilon -\epsilon.
\end{eqnarray*}
It suffices to take $\epsilon = 1/6$ for a contradiction.

Therefore, as the distance between any two
$\lambda$s is bounded below by the modulus of the
zero of $\ft{\chi_\Omega}$ that is nearest to the
origin, there exists a natural number $P$ so that every translate
of $M$ contains at most $P$ elements of $\Lambda$ and, hence,
there exists a natural number $Q$ (we may take $Q = 2NP$)
so that every translate of
$$
\RR\xi+[0,{\ell\over N}]^d
$$
contains at most $Q$ elements of $\Lambda$.

It follows that $\Lambda$ cannot have positive density,
a contradiction as any spectrum of $\Omega$ (which has
volume $1$) must have density equal to $1$.
\Qed

\noindent
{\sc\small
Department of Mathematics, University of Crete, Knossos Ave.,
714 09 Iraklio, Greece.\\
E-mail: {\tt mk@fourier.math.uoc.gr}, {\tt papadim@math.uoc.gr}
}

\end{document}